\documentclass[a4paper,12pt]{article}

\usepackage{amsmath}
\allowdisplaybreaks
\usepackage{amssymb}
\usepackage[authoryear,round]{natbib}
\usepackage{graphicx}
\usepackage[OT4]{fontenc}
\usepackage[utf8]{inputenc}
\usepackage{url}
\usepackage[unicode,debug,colorlinks=true,linkcolor=black,citecolor=black,urlcolor=black,pdfusetitle=true,breaklinks=true]{hyperref}
\newcommand{\abs}[1]{\left| #1 \right|}

\newcommand{\okra}[1]{\left( #1 \right)}

\newcommand{\ceil}[1]{\left\lceil #1 \right\rceil}
\newcommand{\floor}[1]{\left\lfloor #1 \right\rfloor}

\newcommand{\klam}[1]{\left\{ #1 \right\}}

\newcommand{\boole}[1]{{\bf 1}{\klam{#1}}}
\DeclareMathOperator{\card}{\#}
\DeclareMathOperator{\mean}{\mathbf{E}}
\DeclareMathOperator{\var}{Var}

\newcommand{\strong}[1]{{\rm\textbf{#1}}}

\newtheorem{definition}{Definition}
\newtheorem{example}{Example}
\newtheorem{proposition}{Proposition}

\newenvironment*{proof*}[1]{\begin{trivlist}\item[]
\noindent\strong{Proof of #1:}}{$\Box$\par\end{trivlist}}

\title{Vocabulary Growth Fundamentals: \\ 
  Bernstein Functions and Hausdorff Sequences}

\author{{\L}ukasz D\k{e}bowski%
  \\[1ex]
  \small{ORCID: 0000-0001-7136-5283}
  \\
  \small{Institute of Computer
    Science, Polish Academy of Sciences}
  \\
  \small{ul.\ Jana Kazimierza 5,
    01-248 Warszawa, Poland}
  \\
  \small{e-mail: ldebowsk@ipipan.waw.pl}}
\date{}

\begin{document}

\begin{titlepage}
\thispagestyle{empty}
  
\maketitle

\begin{abstract}
  We survey the theory of vocabulary growth founded in the setting of
  stochastic processes. In particular, we model the expected number of
  types through Bernstein functions and Hausdorff sequences. These
  classes of mathematical objects, defined by alternating signs of
  their derivatives or differences, can be related to continuous-time
  Poisson point processes and discrete-time IID processes,
  respectively. Building on previous accounts of the vocabulary
  growth, we integrate the broader theories of Bernstein functions and
  Hausdorff sequences and connect them with recently developed hapax
  rate models. In particular, we prove that the logistic hapax rate
  model has a non-negative spectrum and hence it defines a Bernstein
  function, thereby solving an earlier posed problem. We also analyze
  the limitations of the Bernstein--Hausdorff theory of the vocabulary
  growth by considering its generalizations under stationary and
  Weibull renewal processes.
  \\[2ex]
  \strong{Keywords:} Poisson processes, IID processes, stationary
  processes, Weibull distribution, Bernstein functions, Hausdorff
  moment theorem, hapax rate, Heaps' law, Zipf's law
  \\[1ex]
  \strong{MSC2020:} 60G55 (Primary) 60G10, 91F20 (Secondary)
\end{abstract}
  
\end{titlepage}

\section{Introduction}
\label{secIntroduction}


Two developments may propel quantitative linguistic research in the
coming years. On the one hand, availability of trillion-token corpora
provides an unprecedented empirical ground on which theoretical models
can be tested. On the other hand, mathematical theories such as those
to be presented in this paper, offer multiple models to describe word
frequency distributions. This conjunction of empirical and theoretical
abundance raises a question: What exactly should be modeled in
quantitative linguistics?

As long as we are interested only in marginal distributions, various
sorts of memoryless processes yield a surprisingly good prediction of
language data \citep{Baayen01}. It turns out to be so despite
significant memory effects such as long-range correlations
\citep{Li90, AltmannPierrehumbertMotter09, LinTegmark17,
  TanakaIshii21, MikhaylovskiyChurilov2023, StaniszDrozdzKwapien24,
  BartnickiOther25, WieczynskiDebowski25}, Hilberg's law
\citep{Hilberg90, TakahiraOther16}---related to the neural scaling law
\citep{KaplanOther20, Debowski25g}, and the cubic-logarithmic growth
of the maximal repetition length \citep{Debowski15f, Debowski25e}.

Traditionally, marginal word frequency models have been evaluated in
terms of their ability to support \strong{interpolation}: predicting
marginal frequency distributions of text sections given the statistics
of the whole text \citep{Baayen01, Milicka09, Milicka13, Davis18}.  A
more ambitious and controversial goal has been \strong{extrapolation}:
inferring structural regularities beyond the observed data, as in many
discussions of Zipf's law \citep{Estoup16, Condon28, Zipf35, Zipf49}.
Yet a third perspective, increasingly relevant in the era of massive
corpora, is \strong{mixture identification}: decomposing the lexicon
into statistically distinct components.

Let us work out these three approaches on an example. Let $v(n)$
denote the expected number of types that occur in the first $n$
tokens. For many natural texts one observes an approximate power-law
growth,
\begin{align}
  v(n) \propto n^\beta, \qquad 0<\beta<1,
\end{align}
known as Heaps' law or Herdan's law \citep{Herdan64,Heaps78}. The
first approach, interpolation means predicting $v(m)$ for $m<n$ when
the distribution for $n$ tokens is known. The second approach,
extrapolation is attempting to predict $v(m)$ for $m>n$. The third
approach, mixture identification assumes that the lexicon is a
patchwork of components. Some are frequent, others rare. Formally, one
may write
\begin{align}
v(n) = v_1(n) + v_2(n) + \cdots + v_k(n),
\end{align}
where each $v_i(n)$ corresponds to a particular lexicon component with
a distinct statistical behavior. The goal is to discern those
components.

The question of whether the vocabulary of a natural language is closed
has long been debated. In reality, the lexicon is an evolving entity:
new words are continually coined, while others fall into disuse.  Thus
the lexicon may be closed but only within a moving window. This sounds
like a tautology unless we assume a certain well-defined process of
latent lexicon evolution. A closed vocabulary assumption approximates
certain stable parts of the lexicon but cannot capture the whole. A
more realistic view imagines the lexicon as a mixture of heterogeneous
components: function words, content words, proper names, technical
terms, loanwords, and accidental noise such as typos.  The real
challenge, following Francis Bacon's metaphor of finding ``joints in
nature'', is to determine whether these components indeed behave in
statistically different ways, and if so, where the natural boundaries
lie.

This perspective also reshapes how we interpret diagnostic plots. For
short texts, word frequency distributions and their associated summary
statistics are relatively simply behaved. Let us define the hapax rate
\begin{align}
h(u) = \frac{v(e^u|1)}{v(e^u)},
\end{align}
where $v(n|1)$ is the number of types seen exactly once in a sample of
length $n$ and $u=\log n$.  For novel-sized texts in English, the
hapax rate decreases approximately linearly as a function of variable
$u$ \citep{Debowski25} and the plot of $v(n)$ displays a log-log
concavity \citep{FontClosCorral15}.

However, once we scale to extreme text collections, more complex
patterns emerge. In particular, Zipf's law often splits into two
regimes \citep{FerrerSole01b, MontemurroZanette02, PetersenOther12,
  GerlachAltmann13, WilliamsOther15}. Consequently, the hapax rate
becomes $U$-shaped: $h(u)$ decreases for small $u$ to resurge for
large $u$ \citep{Fan10}. This phenomenon is consistent with mixture
models \citep{Debowski25} that combine a monotone decreasing hapax
rate component, likely implying a closed lexicon, with a nuisance
component, possibly power-law distributed as in the monkey-typing
explanation of Zipf's law \citep{Mandelbrot54, Miller57}. Such
mixtures reflect the intuition that in massive corpora, besides the
systematic structure of natural language, there also accumulate rare
types that many linguists would not call ``words'' at all.


The goal of the paper is to supply a mathematical framework that can
address challenges of linguistic big data. This framework is founded
on Hausdorff sequences and Bernstein functions, which turn out to be
generic models of the expected number of types $v(n)$. Abstractly
speaking, Hausdorff sequences and Bernstein functions are arbitrary
non-negative sequences or smooth functions whose subsequent
differences or derivatives have alternating signs \citep{Hausdorff23,
  Bernstein28, SchillingSongVondracek10, Akhiezer21}. This concept
arises naturally if we consider discrete-time IID processes and
continuous-time Poisson point processes \citep{Kingman67,
  Kingman93}. In both cases, intensities of word types are distributed
according to a non-negative measure given by the Hausdorff moment
theorem \citep{Hausdorff23} or by the L\'evy-Khintchine representation
\citep{Bernstein28, Levy48}.  In particular, \citet{Baayen01}
presented an introduction to such models for quantitative linguists,
preceded by more mathematical works by \citet{Karlin67} and
\citet{Khmaladze88}.

Our aim is to extend the discussion by \citet{Baayen01} by integrating
that framework with a more general theory of Hausdorff sequences
\citep{Akhiezer21} and Bernstein functions
\citep{SchillingSongVondracek10}. Also, we take into account recent
research in parametric hapax rate models \citep{Debowski25}. The
latter study offered new perspectives on modeling marginal word
frequency distributions but was unaware of more powerful results.
Besides presenting connections that are known, we prove several new
results. In particular, we solve a problem stated by
\citet{Debowski25}. Namely, we prove that the logistic hapax rate
model, which gives a good prediction for novel-sized texts, has
non-negative spectrum elements and, thus, it yields a Bernstein
function. Another new item is a stress test of the
Bernstein--Hausdorff theory of the vocabulary growth. Namely, we
investigate what survives of this theory when we replace memoryless
sources by a stationary process or a Weibull-distributed renewal
process.

This paper is a theoretical survey. We disclose that it was inspired
by a dialogue with ChatGPT about relevant mathematical concepts. That
is, the chatbot was used as an intelligent search engine for abstract
patterns appearing in the training data---existing mathematical
literature. Intentionally, we do not analyze empirical data in this
paper but we rather review mathematical theory that is scattered over
various sources and essentially known since a long time.  Our aim is
to stimulate further cross-disciplinary dialogue and encourage other
language researchers to read out of their domain and apply the results
mentioned here. The amount of theory to catch up with is enough as for
a full paper.

The organization of the paper is as follows.  Section \ref{secSetting}
introduces four kinds of stochastic processes that generate word
tokens: the Poisson source analyzed in Section \ref{secPoisson}, the
multinomial source analyzed in Section \ref{secMultinomial}, the
invariant source treated in Section \ref{secInvariant}, and the
Weibull source treated in Section \ref{secWeibull}. The paper is
concluded in Section \ref{secConclusion}. Appendices
\ref{secPoissonA}--\ref{secWeilbullA} contain the proofs of
propositions stated in Sections \ref{secPoisson}--\ref{secWeibull},
respectively.

\section{Stochastic setting}
\label{secSetting}


Let us take a bird's eye perspective onto a spoken text in natural
language.  At different time positions in the text, which can be
expressed by real numbers, we encounter various word tokens. We are
able to group these tokens into certain types
$w\in\mathbb{W}$---typically representing distinct words. In
quantitative linguistics, we are interested in counting types, i.e.,
the number of distinct words, and the spectrum, i.e., the numbers of
word types that occur $k$ times in the text, where
$k=1,2,\ldots$. Before we specify a particular stochastic source $S$
that models the distribution of word tokens, we will introduce some
general notations.

\subsection{Counting types}
\label{secCounting}

Let $N^S_w((a,b])$ be the random number of occurrences of word type
$w$ for a process $S$ in the time interval $(a,b]$. Function
$N^S_w((a,b])$ is an instance of a random \strong{counting measure}
\citep{Kingman93}. In particular,
$N^S_w((a,c])=N^S_w((a,b])+N^S_w((b,c])$.  We define:
\begin{itemize}
\item $N_S(t):=\sum_{w\in\mathbb{W}} N^S_w((0,t])$ ---
  the \strong{number of tokens},
\item
  $V_S(t):=\sum_{w\in\mathbb{W}} \boole{N^S_w((0,t])>0}$
  --- the \strong{number of types},
\item
  $V_S(t|k):=\sum_{w\in\mathbb{W}}
  \boole{N^S_w((0,t])=k}$ --- the \strong{spectrum element}.
\end{itemize}
We observe some simple constraints such as $V_S(t|k)\ge 0$ and
\begin{align}
  \label{Consistent}
  \sum_{k=1}^\infty V_S(t|k)&=V_S(t),
  &
    \sum_{k=1}^\infty kV_S(t|k)&=N_S(t). 
\end{align}
In particular, $V_S(t|1)$ is the number of types that occur once. These
are called \emph{hapax legomena} in Greek or, succinctly, hapaxes in
English.

From the above entities, we can compute the rank-frequency plot. The
maximal rank with frequency $f$, called the \strong{rank function},
can be expressed and bounded  as
\begin{align}
  R_S(t|f)
  &:=\sum_{w\in\mathbb{W}}
    \boole{N^S_w((0,t])\ge f}
    =
    V_S(t)-\sum_{k=1}^{f-1} V_S(t|k)
    \le
    \frac{N_S(t)}{f},
\end{align}
see the work of \citet{Debowski25}. The harmonic bound on the right
reminds of Zipf's law \citep{Zipf35} and stems from the Markov
inequality.

Among these objects, spectrum elements play the central role in
modeling word frequency distributions, as all other entities can be
computed from it.  We denote the \strong{expectation} and the
\strong{variance}
\begin{align}
  \mean X&:=\int XdP, & \var X&:=\mean (X-\mean X)^2.
\end{align}
As we will see further, the expected spectrum can be
also computed from the expected number of types exactly or
approximately, depending on a particular source.

\subsection{Four kinds of sources}
\label{secFourSources}

There are four basic idealized sources of word tokens. These are:
\begin{itemize}
\item the \strong{invariant source} (I): a general stationary sequence
  of words;
\item the \strong{multinomial source} (M): an IID sequence of words;
\item the \strong{Weibull source} (W): a collection of independent IID
  sequences of interword intervals, being a special case of the
  renewal process with the Weibull distribution;
\item the \strong{Poisson source} (P): a particularly symmetric special
  case of the Weibull source that reduces to a collection of Poisson
  point processes.
\end{itemize}
Each of the above four processes induces a specific collection of
\strong{counting measures} $(N^S_w)_{w\in\mathbb{W}}$, where
$S\in\klam{I,M,W,P}$. Processes (M) and (P) were considered by
\citet{Baayen01}.  Processes (M) and (I) generate word tokens at
integer positions, one word per position. By contrast, processes (W)
and (P) generate word tokens at real positions, so that any open
interval may contain none, one, or multiple words. Let us inspect
these four sources.

To begin,  a \strong{stationary process} is a
sequence $(X_i)_{i\in\mathbb{Z}}$ of random variables such that any
tuple $(X_{i+1},\ldots,X_{i+k})$ has the same distribution as
$(X_{1},\ldots,X_{k})$.  The invariant source is a stationary process
$(X_i)_{i\in\mathbb{Z}}$ with variables taking values in word types,
$X_i:\Omega\mapsto\mathbb{W}$. For the stationary process
$(X_i)_{i\in\mathbb{Z}}$, we define the \strong{invariant counting
  measures} $(N^I_w)_{w\in\mathbb{W}}$ on the real line as
\begin{align}
  \label{Counting}
  N^I_w((a,b]):=\card\klam{i\in (a,b]:X_i=w},
  \quad a,b\in\mathbb{R}.
\end{align}

Some special case of a stationary process is an \strong{IID process},
which is a sequence $(Z_i)_{i\in\mathbb{Z}}$ of independent
identically distributed (IID) random variables. Consider an IID
process $(X_i)_{i\in\mathbb{Z}}$ with the marginal distribution
$P(X_i=w)=p_w$, where
\begin{align}
  \label{Probability}
  p_w\ge 0, \quad \sum_{w\in\mathbb{W}} p_w=1.
\end{align}
In this case, the invariant counting measures
$(N^I_w)_{w\in\mathbb{W}}$ reduce to the \strong{multinomial counting
  measures} $(N^M_w)_{w\in\mathbb{W}}$. The distribution of the random
variable $N^M_w((0,t])$ is given by the \strong{binomial distribution}
\begin{align}
  \label{Binomial}
  P(N^M_w((0,t])=k)=
  \binom{\floor{t}}{k} p_w^k(1-p_w)^{\floor{t}-k},
  \quad t\ge 0.
\end{align}

The assumption that word tokens are probabilistically independent is,
of course, a crude approximation. Some simple adjustment is to assume
that the length of intervals between the occurrences of the same word
type are an IID process. This yields a renewal process. Typically in
language modeling, it is assumed that successive recurrence times
$(Y^w_i)_{i\in\mathbb{Z}}$ of a word type $w$ are a real-valued IID
process with the marginal \strong{Weibull distribution}
\begin{align}
  \label{Weibull}
  P(Y^w_i\le y)=
  \int_0^y
  \frac{\alpha_w}{\lambda_w}\okra{\frac{t}{\lambda_w}}^{\alpha_w-1}
  e^{-(t/\lambda_w)^{\alpha_w}}dt.
\end{align}
Parameter $\lambda_w$ is called the scale, whereas parameter
$\alpha_w$ is called the shape---or the burstiness in quantitative
linguistics. Burstiness describes a tendency of a word type to occur
in clusters.  Typically $\alpha_w<1$ for words in natural language and
there is a huge variation of the estimates of $\alpha_w$ across
different word types \citep{AltmannPierrehumbertMotter09}.

For a collection of independent IID processes
$(Y^w_i)_{i\in\mathbb{Z}}$ distributed along (\ref{Weibull}) for
$w\in\mathbb{W}$, we define the \strong{Weibull counting measures}
$(N^W_w)_{w\in\mathbb{W}}$ on the real line as
\begin{align}
  N^W_w((0,t]):=\card\klam{k\ge 0:\sum_{i=1}^k Y_i\le t},
  \quad t> 0.
\end{align}

The Weibull source has a particularly symmetric case.  For
$\alpha_w=1$ and $\lambda_w=1/p_w$ where intensities $p_w$ satisfy
condition (\ref{Probability}), the Weibull counting measures
$(N^W_w)_{w\in\mathbb{W}}$ reduce to the \strong{Poisson counting
  measures} $(N^P_w)_{w\in\mathbb{W}}$. Equivalently, Poisson counting
measures are induced by a collection of probabilistically independent
Poisson point processes \citep{Kingman93} and therefore, they enjoy a
few nice properties. In particular, variables $N^P_w((a,b])$ and
$N^P_w((c,d])$ with $a,b,c,d\in\mathbb{R}$ are probabilistically
independent for $(a,b]\cap (c,d]=\emptyset$ and variable
$N^P_w((s,s+t])$ has the \strong{Poisson distribution}
\begin{align}
  \label{Poisson}
  P(N^P_w((s,s+t])=k)=\frac{(p_w t)^k}{k!} e^{-p_w t},
  \quad t> 0.
\end{align}
For large interval lengths $t$ and small probabilities $p_w$, the
Poisson distribution (\ref{Poisson}) approximates the binomial
distribution (\ref{Binomial}). Hence the Poisson counting measures are
a certain approximation of the multinomial counting measures that may
be good for long samples.

\section{Poisson sources}
\label{secPoisson}

We begin the analysis with the Poisson source since it is the most
symmetric case and enjoys the most developed theory. The expectations
of our quantities of interest can be computed exactly.
\begin{proposition}
  \label{theoMean}
  For the Poisson source, under condition (\ref{Probability}) and $t>0$, we
  have
  \begin{align}
    \label{Tokens}
    \mean N_P(t)
    &=t,
    \\
    \label{VarTokens}
    \var N_P(t)
    &=t,
    \\
    \label{SmallV}
    \mean V_P(t)
    &=v_P(t)
    := \sum_{w\in\mathbb{W}} (1-e^{-p_w t}),
    \\
    \label{VarV}
    \var V_P(t)
    &=v_P(2t)-v_P(t)\le v_P(t),
    \\
    \label{SmallVK}
    \mean V_P(t|k)
    &=v_P(t|k)
      := \sum_{w\in\mathbb{W}} \frac{(p_w t)^k}{k!} e^{-p_w t},
    \\
    \label{VarVK}
    \var V_P(t|k)
    &=v_P(t|k)-\frac{1}{2^{2k}}\binom{2k}{k}v_P(2t|2k)\le v_P(t|k).
  \end{align}
\end{proposition}
In formula (\ref{VarVK}), we may bound
$\frac{1}{2k+1}\le\frac{1}{2^{2k}}\binom{2k}{k}\le 1$.

Now we will see that the expected spectrum can be computed from the
expected number of types by taking derivatives. In the following, we
write the \strong{derivative operator}
\begin{align}
  \label{Derivatives}
  \nabla f(t)&:=\frac{df(t)}{dt}.
\end{align}
We also write $\nabla^0 f(t):=f(t)$ and
$\nabla^{k+1} f(t):=\nabla \nabla^k f(t)$.  We assume that this
operator acts only on the directly following functor:
$\nabla f(g(t))=df(g(t))/dg(t)$.
\begin{proposition}
  \label{theoTaylor}
  For the Poisson source, under $n\in\mathbb{N}$ and
  $t>0$, we have
  \begin{align}
    \label{Taylor}
    v_P(t|k)
    &=(-1)^{k+1}\frac{t^k}{k!}\nabla^k v_P(t), \qquad k\in\mathbb{N}.
  \end{align}
\end{proposition}

\subsection{Bernstein functions}
\label{secBernsteinI}

Successful modeling of empirical word frequency distributions rests on
supplying some relatively simple models of function $v_P(t)$ while the
exact intensities of word types remain unknown.  For arbitrary
function $v_P(t)$, let us treat claim (\ref{Taylor}) as the definition
of functions $v_P(t|k)$. Consequently, conditions $v_P(t)\ge 0$ and
$v_P(t|k)\ge 0$ imply that functions and $v_P(t)$ is non-negative and
its subsequent derivatives have alternating signs. Such a concept is
well known in mathematics and gives rise to a profound theory with
multiple applications in engineering and statistics.
\begin{definition}[\mbox{\citealp[Def.\ 1.3 and Def.\
    3.1]{SchillingSongVondracek10}}]
  A smooth function $v:(0,\infty)\to\mathbb{R}$ is called a
  \strong{Bernstein function} if $v(t)\ge 0$ and
  $(-1)^{k+1}\nabla^k v(t)\ge 0$ for all $k\in\mathbb{N}$ and $t>0$.
  A smooth function $u:(0,\infty)\to\mathbb{R}$ is called
  \strong{completely monotone} if $u(t)=\nabla v(t)$ for a certain
  Bernstein function $v:(0,\infty)\to\mathbb{R}$ and all $t>0$.
\end{definition}
For any function $v:(0,\infty)\to\mathbb{R}$, we write
$v(0+):=\lim_{t\to 0+} v(t)$. We call a Bernstein function $v(t)$
\strong{standard} if $v(0+)=0$ and $\nabla v(0+)=1$. A Bernstein function
$v(t)$ with $\nabla v(0+)=\infty$ is called \strong{singular}.

Let us present some simple examples of Bernstein functions.
\begin{example}
  Obviously, function $v_P(t)$ defined in (\ref{SmallV}) is a
  Bernstein function.  Here are some more examples of Bernstein
  functions: $t\mapsto t^\beta$,
  $t\mapsto \beta^{-1}\okra{(t+1)^\beta-1}$, $t\mapsto\log(t+1)$, and
  $t\mapsto p^{-1}\okra{1-e^{-p t}}$ for $\beta\in(0,1]$ and $p>0$.
  The first one is singular, whereas the other three are standard.
\end{example}

\subsection{Bernstein theorem}
\label{secBernstein}

It is true that any Bernstein function is a mixture of elementary
functions $t\mapsto p^{-1}\okra{1-e^{-p t}}$ for varying $p>0$. A
particular case of a discrete mixture has been already encountered in
formula (\ref{SmallV}).  The more general result is known as the
Bernstein theorem on completely monotone functions \citep{Bernstein28}
or the L\'evy-Khintchine representation \citep{Levy48}.  In the
following statement of this result, we apply the \strong{Lebesgue
  integral} \citep{Billingsley79}.
\begin{proposition}[\mbox{\citealp[Thm.\ 1.4 and Thm.\
    3.2]{SchillingSongVondracek10}}]\hfill
  \label{theoVF}
  A smooth function $v:(0,\infty)\to\mathbb{R}$ is a Bernstein
  function if and only if there exists a unique non-negative measure
  $F$ on $[0,\infty)$ such that
  $F(\klam{0})+\int_{(0,\infty)}\min\klam{1,p^{-1}}dF(p)<\infty$ and
  \begin{align}
    \label{FtoV}
    v(t)&=v(0+)
          +
          tF(\klam{0})
          +
          \int_{(0,\infty)}\frac{1-e^{-p t}}{p}dF(p),
    \\
    \label{FtoVPrime}
    \nabla v(t)&=F(\klam{0})
          +
          \int_{(0,\infty)}e^{-p t}dF(p).
  \end{align}
\end{proposition}
Additionally, we see that a Bernstein function
$v:(0,\infty)\to\mathbb{R}$ is standard if and only if $v(0+)=0$ and
the measure $F$ is a \strong{probability measure}, namely,
\begin{align}
  F(\klam{0})+\int_{(0,\infty)} dF(p)=1.
\end{align}
Formula (\ref{FtoVPrime}) states that the derivative of $v(t)$ is the
\strong{Laplace transform} of measure $F$. In probability theory,
measure $G$ given by $dG(p)=dF(p)/p$ is called the L\'evy measure of
function $v(t)$. In the linguistic context, measure $G$ has been
called the structural type distribution, whereas $F$ has been called
the structural token distribution by \citet{Baayen01}. Preferring
short names, we will refer to $F$ as the \strong{Bernstein measure}.

\begin{example}
  The Bernstein function $v(t)=v_P(t)$ defined in (\ref{SmallV}) has a
  Bernstein measure (\ref{Points}) that is supported on discrete points.
\end{example}

Many useful Bernstein functions enjoy a continuous Bernstein measure that
has a density $f:(0,\infty)\to\mathbb{R}$ with respect to the Lebesgue
measure.  We call $f$ the \strong{Bernstein density}.  Then we may
substitute
\begin{align}
  dF(p)=f(p)dp.
\end{align}
\begin{example}
  In general, the Bernstein measure can be supported
  out of the unit interval $[0,1]$, which corresponds to occurrence of
  ultra-frequent types with intensity $p_w>1$. For example, the
  Bernstein function $v(t)\propto t^\beta$, where $\beta\in(0,1)$, is
  yielded by the Bernstein density $f(p)\propto p^{-\beta}$.
\end{example}
We will inspect more examples in Section \ref{secModels}.

\subsection{Taylor expansion}
\label{secTaylor}

Let us make a short excursion to complex analysis.  It is known that
the Laplace transform is often an analytic function on the positive
complex half-plane, which can be a useful fact. In the case of
Bernstein functions we have a somewhat stronger statement.
\begin{proposition}[\mbox{\citealp[Proposition
    3.5]{SchillingSongVondracek10}}]
  \label{theoAnalyticV}
  Any Bernstein function $v:(0,\infty)\to\mathbb{R}$ has an extension
  $v(z)$ that is analytic for $\Re z>0$ and continuous for
  $\Re z\ge 0$.
\end{proposition}

This result can be applied to the spectrum elements defined as in
(\ref{Taylor}).
\begin{definition}
  For a smooth function $v:(0,\infty)\to\mathbb{R}$, we define the
  \strong{spectrum elements}
  \begin{align}
  \label{TaylorII}
  v(t|k)&:=(-1)^{k+1}\frac{t^k}{k!}\nabla^k v(t), & k\in\mathbb{N}.
  \end{align}
\end{definition}
We see that the spectrum elements are the Taylor expansion
coefficients of the Bernstein function. Since the Bernstein function
is analytic in the positive half-plane, the respective Taylor series
converges to the function on the maximal disk contained in the
positive half-plane. This gives rise to a useful representation.
\begin{proposition}
  \label{theoPowerSeries}
  For the extension $v(z)$ of a Bernstein function $v(t)$ and any
  $t>0$, we have
  \begin{align}
    \label{PowerSeries}
    v(t-z t)
    &=
      v(t)-\sum_{k=1}^\infty z^k v(t|k), \quad \abs{z}\le 1,
    \\
    \label{PowerSeriesPrime}
    t\nabla v(t-z t)
    &=
      \sum_{k=1}^\infty k z^{k-1} v(t|k), \quad \abs{z}<1 \lor z=1.
  \end{align}
\end{proposition}

Applying formula (\ref{PowerSeries}) for $z=\pm 1$ to formulas
(\ref{VarV}), we obtain symmetric expressions for the expectation and
the variance of the number of types in the marked Poisson case
\begin{align}
  \label{PoissonVMean}
  \mean V_P(t)=v_P(t)-v_P(0+)&=\sum_{k=1}^\infty (+1)^{k-1} v_P(t|k),
  \\
  \label{PoissonVVar}
  \var V_P(t)=v_P(2 t)-v_P(t)&=\sum_{k=1}^\infty (-1)^{k-1} v_P(t|k).
\end{align}

By formulas (\ref{PowerSeries}) and (\ref{PowerSeriesPrime}),
consistency conditions analogous to
(\ref{Consistent}) can be asserted easily for general Bernstein
functions.
\begin{proposition}
  \label{theoTypeConsistent}
  A Bernstein function $v:(0,\infty)\to\mathbb{R}$ satisfies condition
  $\sum_{k=1}^\infty v(t|k)=v(t)$ if and only if $v(0+)=0$.
\end{proposition}
\begin{proposition}
  \label{theoTokenConsistent}
  A Bernstein function $v:(0,\infty)\to\mathbb{R}$ satisfies condition
  $\sum_{k=1}^\infty kv(t|k)=t$ if and only if $\nabla v(0+)=1$.
\end{proposition}
Thus only standard Bernstein functions satisfy both consistency
conditions and can be applied as models of word frequency
distributions.

\subsection{Hapax rate}
\label{secHapax}

We adopt this definition.
\begin{definition}
  For a Bernstein function $v:(0,\infty)\to\mathbb{R}$, we define
  the \strong{hapax rate} $h:\mathbb{R}\to[0,1]$ where
  \begin{align}
    \label{Rate}
    h(u)&:=\frac{v(e^u|1)}{v(e^u)-v(0+)}.
  \end{align}
\end{definition}

According to the next statement, there is a one-to-one correspondence
between hapax rates and Bernstein functions if we fix the value of the
latter at lengths $t\in\klam{0,1}$. Formula (\ref{HtoV}) follows by
definition (\ref{TaylorII}) and solving a simple ordinary differential
equation.
\begin{proposition}[\mbox{\citealp[\S 3.4]{Debowski25}, \citealp[Eq.\
    (3.43)]{Baayen01}}]
  \label{theoVH}
  For a smooth function $h:\mathbb{R}\to[0,1]$ and a smooth function
  $v:(0,\infty)\to\mathbb{R}$, condition (\ref{Rate}) holds if and
  only if
  \begin{align}
    \label{HtoV}
    v(t)-v(0+)&=(v(1)-v(0+))\exp\okra{\int_0^{\log t} h(u)du},
    \\
    \label{HKtoVK}
    v(t|k)&=(v(t)-v(0+))h(\log t|k),
  \end{align}
  where
  \begin{align}
    \label{HK1}
    h(u|0)&:=-1,
    \\
    \label{HK2}
    h(u|k)&:=\okra{1-\frac{1}{k}\okra{1+h(u)+\frac{d}{du}}}h(u|k-1),
              \quad k\ge 1.
  \end{align}
\end{proposition}
Functions $h(u|k)$ are called the \strong{relative spectrum elements}
\citep[p.\ 90]{Baayen01}. Formula (\ref{HtoV}) was observed by
\citep[\S 3.4]{Debowski25}. Recursion (\ref{HK2}) was observed by
\citet[Eq.\ (3.43)]{Baayen01}.

Consequently, we can characterize hapax rates in this way.
\begin{proposition}[\mbox{\citealp[\S 3.4]{Debowski25}}]
  \label{theoH}
  A smooth function $h:\mathbb{R}\to[0,1]$ is the hapax rate of a
  Bernstein function if and only if
  \begin{align}
    \label{RetroH}
    \int_{-\infty}^0 h(u)du&=\infty,
    \\
    \label{PositiveHK}
    h(u|k)&\ge 0, \quad k\ge 1.
  \end{align}
\end{proposition}
\begin{proposition}
  \label{theoHStandard}
  Under conditions (\ref{RetroH}) and
  (\ref{PositiveHK}), function $v(t)$ given by (\ref{HtoV}) defines a
  standard Bernstein function if and only if
  \begin{align}
    \label{RetroDampedH}
    \int_{-\infty}^0 (1-h(u))du&<\infty,
    \\
    \label{V0}
    v(0+)&=0,
    \\
    \label{V1}
    v(1)&=\exp\okra{-\int_{-\infty}^0(1-h(u))du}.
  \end{align}
\end{proposition}
Thus, there is a one-to-one correspondence between hapax rates that
satisfy condition (\ref{RetroDampedH}) and standard Bernstein
functions. Condition (\ref{RetroDampedH}) holds if and only if
$\nabla v(0+)<\infty$.

Asserting property (\ref{PositiveHK}) for a particular Bernstein
function is a task that requires some algebraic invention. We note
that conditions $h(u|1)\ge 0$ and $h(u|2)\ge 0$ are equivalent to
\begin{align}
  \label{HHPrime}
  h(u)&\ge 0,
  & \frac{dh(u)}{du}&\le h(u)(1-h(u)).
\end{align}
\begin{example}
  For the logistic hapax rate $h(u)=1/(e^{\gamma u}+1)$ with a
  parameter $\gamma\in(0,1]$, to be discussed in Section
  \ref{secModels}, we obtain $dh(u)/du=-\gamma h(u)(1-h(u))$ and the
  relative spectrum elements $h(u|k)$ can be expressed as polynomials
  of variable $h(u)$ \citep{Debowski25}.
\end{example}

\subsection{Inverse Laplace transform}
\label{secInverse}

Finding the Bernstein measure for a given Bernstein function $v(t)$
can be hard since measure $F$ given by formulas
(\ref{FtoV})--(\ref{FtoVPrime}) is the inverse Laplace transform of
the derivative $\nabla v(t)$. Besides using look-up tables of Laplace
transforms, as the last resort, we can use the Post-Widder inversion
formula \citep{Post30, Widder41}, which works if the Bernstein density
$f(p)$ at a point $p$ exists.  Applying the spectrum elements
$v(t|k)$, the formula reads simply
\begin{align}
  f(p)
  &=
    \lim_{k\to\infty} k v\okra{\frac{k}{p}\middle|k}
    \nonumber\\
  &=
  \lim_{k\to\infty} k v\okra{\frac{k}{p}}
  h\okra{\log k-\log p\middle|k},
  \label{PostWidder}
\end{align}
so $f(p)\ge 0$ if (\ref{PositiveHK}) holds. This can be used for
instance to infer $f(p)\propto p^{-\beta}$ for $v(t)\propto t^\beta$
and $\beta\in(0,1)$ but it is far less effective for other examples.

The Post-Widder inversion formula (\ref{PostWidder}) cannot be used to
infer the value $F(\klam{0})$. However, we can recover 
\begin{align}
  F(\klam{0})=\lim_{t\to\infty} \frac{v(t)}{t},
\end{align}
as commented by \citet[p.\ 16, Remark
(iv)]{SchillingSongVondracek10}.  Rate $v(t)/t$ is called the
type-token ratio (TTR) in quantitative linguistics. By formulas
(\ref{HtoV}) and (\ref{V1}), for a standard Bernstein function we
obtain
\begin{align}
  \label{TTR}
  \frac{v(t)}{t}&=\exp\okra{-\int_{-\infty}^{\log t}(1-h(u))du}.
\end{align}
Formula (\ref{TTR}) predicts that the type-toke ratio is monotonically
decreasing with the text length. The further discussed Proposition
\ref{theoBlockTypes} generalizes this property to stationary sources.

\subsection{Bounded and expansive functions}
\label{secBoundedExpansive}

Under the Bernstein function framework, a lexicon component that is
closed corresponds to the notion of a bounded Bernstein function.
\begin{definition}
  A Bernstein function $v(t)$ is called \strong{bounded} if
  \begin{align}
    \label{ClosedV}
    \lim_{t\to \infty} v(t)<\infty.
  \end{align}
\end{definition}
By formula (\ref{HtoV}), a Bernstein function $v(t)$ is bounded if and
only if the corresponding hapax rate $h(u)$ satisfies condition
\begin{align}
  \label{ClosedH}
  \int_0^\infty h(u)du&<\infty.
\end{align}
So, the lexicon component is closed when the area under the hapax rate
plot is finite.

By contrast, a lexicon component that keeps growing across all
spectrum elements corresponds to a concept of an expansive Bernstein
function.
\begin{definition}
  A Bernstein function $v(t)$ is called \strong{expansive} at a point
  $t>0$ if
  \begin{align}
    \label{GrowingVK}
    \nabla v(t|k)&\ge 0, \quad k\ge 1.
  \end{align}
\end{definition}

Expression (\ref{PoissonVVar}) can be also used to
prove a fact that sounds obvious but is not immediate. Namely,
expansive functions cannot be bounded.
\begin{proposition}
  \label{theoBoundedExpansive}
  A Bernstein function cannot be both bounded and expansive at all
  points.
\end{proposition}

\subsection{Preserving Bernstein property}
\label{secTransformations}

Much work in the theory of Bernstein functions revolves around
transformations that preserve their defining property
\citep{SchillingSongVondracek10}. There are many of them and we will
only name a few basic ones.

The first transformation is a linear mapping of the $x$ and $y$ axes.
\begin{definition}
  Consider parameters $a,b,c\ge 0$ and $d\ge -cv(b)$.  The
  \strong{affine transformation} of a function $v(t)$ is defined as
  function
  \begin{align}
    v_{abcd}(t)=cv(a t+b)+d.
  \end{align}
\end{definition}
The affine transformation of a Bernstein function is a Bernstein
function.  In particular, the affine transformation is standard if we
put
\begin{align}
  c&=\frac{1}{a\nabla v(b)}, & d&=-\frac{v(b)}{a\nabla v(b)}.
\end{align}
\begin{example}
  In this way, a non-standard Bernstein function $t\mapsto t^\beta$
  gives rise to a standard Bernstein function
  $t\mapsto \beta^{-1}\okra{(t+1)^\beta-1}$ for $\beta\in(0,1]$.
\end{example}

The next transformation is a weighted average of lexicon components.
\begin{definition}
  Consider a family of functions $v_e(t)$, where $e\in E$. Let $\pi$
  be a probability measure on $E$, namely, $\int_E \sigma(de)=1$. The
  \strong{mixture} of functions $v_e(t)$ is function
  \begin{align}
    \label{Mixture}
    v(t)&=\int_E v_e(t)d\sigma(e).
  \end{align}
\end{definition}
The mixture of Bernstein functions is a Bernstein function. The
mixture is standard if all input functions are standard.

We can also consider limits of Bernstein functions.
\begin{definition}
  Consider a sequence of functions $v_m(t)$, $m\in\mathbb{T}$, that
  converges pointwise for any $t\in(0,\infty)$. The \strong{limit} of
  functions $v_m(t)$ is
  \begin{align}
    v(t)&=\lim_{m\to\infty} v_m(t).
  \end{align}
\end{definition}
The limit of Bernstein functions is a Bernstein function if the
sequence converges uniformly on every interval $[a,b]\in(0,\infty)$
since we can interchange the limit and the differentiation
\citep[Theorem 7.17]{Rudin74}. The limit is standard if all functions
in the sequence are standard.
\begin{example}
  In particular, the limit of functions
  $t\mapsto \beta^{-1}\okra{(t+1)^\beta-1}$ for $\beta\to 0$ is the
  Bernstein function $t\mapsto\log (t+1)$.
\end{example}

Another useful transformation involves composition.  The
\strong{composition} of functions $u(t)$ and $v(t)$ is function
\begin{align}
  u\circ v(t)&=u(v(t)).
\end{align}
As shown by \citet[Theorem 3.6]{SchillingSongVondracek10}, the
composition $u\circ v$ of a completely monotone function $u$ and a
Bernstein function $v$ is completely monotone. Since completely
monotone functions are derivatives of Bernstein functions, for any
$a\ge 0$ and two Bernstein functions $v_1(t)$ and $v_2(t)$, we obtain
another Bernstein function
\begin{align}
  \label{IntegralComposition}
  v(t)=a+\int_0^t \nabla v_1(v_2(m))dm.
\end{align}
Function $v(t)$ is standard if $a=0$, $\nabla v_1(0+)=1$, and $v_2(0+)=0$.
\begin{example}
  \label{exStable}
  In particular, the composition of a completely monotone function
  $t\mapsto e^{-c t}$ and a Bernstein function $t\mapsto t^\gamma$
  yields a completely monotone function $t\mapsto\exp(-c t^\gamma)$ for
  $\gamma\in(0,1]$ and $c>0$, which according to Proposition
  \ref{theoVF}, can be written as
  \begin{align}
    \label{Stable}
    \exp(-c t^\gamma)=\int_0^\infty e^{-p t}g_\gamma(p;c)dp
  \end{align}
  for a certain probability density denoted $g_\gamma(p;c)$, called
  the stable density. This density will be used in the analysis of the
  logistic hapax rate $h(u)=1/(e^{\gamma u}+1)$ in Section
  \ref{secModels}.
\end{example}

\subsection{Mixture identification}
\label{secMixture}

As we wrote in the introduction, an important quantitative linguistic
problem is identifying lexicon components $v_i(t)$ for a given
vocabulary growth function $v(t)$. The corresponding mathematical
problem is as follows: We have a mixture $v(t)$ given by formula
(\ref{Mixture}) for a certain family of Bernstein functions $v_e(t)$, where
$e\in E$. Is the mixture measure $\sigma$ on $E$ given uniquely?

The answer to this question depends of the family of functions
$(v_e(t))_{e\in E}$ but the general method is determined by the
Choquet theorem. The setting is as follows.  We
extend the parameter set from a set $E$ to a convex subset $S$ of a
linear space.  Set $S$ is called \strong{convex} if for any
$s,s'\in S$ and $0\le\lambda\le 1$ we have
$\lambda s+(1-\lambda) s'\in S$.
\begin{definition}
  For a convex set $S$, parametrization $S\ni s\mapsto v_s(t)$ is
  called \strong{affine} if for any $s,s'\in S$ and $0\le\lambda\le 1$
  we have
  \begin{align}
    v_{\lambda s+(1-\lambda) s'}(t)=\lambda v_s(t)+(1-\lambda)v_{s'}(t)
  \end{align}
\end{definition}
For a convex set $S$, point $e\in S$ is called \strong{extremal} if
for any $s,s'\in S$ and $0\le\lambda\le 1$ condition
$\lambda s+(1-\lambda) s'=e$ implies $s=s'=e$. The set of extremal
points is denoted $E(S)$.  The Choquet theorem is as follows.
\begin{proposition}[\mbox{\citealp[p.\ 14]{Phelps01}}]
  \label{theoChoquet}
  For a compact convex set $S$, an affine parametrization
  $S\ni s\mapsto v_s(t)$, and any point $s\in S$, there exists a
  probability measure $\sigma$ on $E(S)$ such that
  \begin{align}
  \label{Choquet}
    v_s(t)=\int_{E(S)} v_e(t)d\sigma(e).
  \end{align}
\end{proposition}

Uniqueness of measure $\sigma$ is a property that we desire for
linguistic applications of decomposing the lexicon into identifiable
parts. In general, uniqueness of measure $\sigma$ is not
guaranteed. For example, consider a ball
$S=\klam{(x,y,z)\in\mathbb{R}^3: x^2+y^2+z^2\le 1}$ in a
three-dimensional space and a sphere
$E(S)=\klam{(x,y,z)\in\mathbb{R}^3: x^2+y^2+z^2=1}$, the boundary of
$S$. The origin $(0,0,0)\in S$ can be expressed as a convex
combination of points on sphere $E(S)$ with a uniform measure $\sigma$
or with a measure $\sigma$ concentrated on two opposite points, say
$(1,0,0)\in E(S)$ and $(-1,0,0)\in E(S)$. By contrast, if we let a
tetrahedron $S=\klam{(x,y,z)\in\mathbb{R}^3: x+y+z\le 1, x,y,z\ge 0}$
then the set of extremal points is
$E(S)=\klam{(0,0,0),(1,0,0),(0,1,0),(0,0,1)}$. For these $S$ and
$E(S)$, any point in $(x,y,z)\in S$ can be written as a unique mixture
(convex combination) of the four points in $E(S)$.

There are some important cases when representation (\ref{Choquet}) is
unique. For example, the Bernstein theorem stated in Proposition
\ref{theoVF} assures such uniqueness for the convex set of Bernstein
functions and extremal points being functions
$t\mapsto p^{-1}\okra{1-e^{-p t}}$ with $p>0$. Another result of a
linguistic importance is the ergodic decomposition for stationary
processes over a countable alphabet, which assures a unique
representation for the convex set of stationary measures and extremal
points being ergodic measures \citep[\S 4.4]{Debowski21}. In general,
representation (\ref{Choquet}) is unique if the compact convex set $S$
is a Choquet simplex---an object that generalizes finite-dimensional
simplices (intervals, triangles, tetrahedra, etc.).

A natural question in the context of quantitative linguistics is
whether mixtures of power laws are also identifiable.  Applying
Proposition \ref{theoVF} and a simple change of variable, we can
easily show that the following is true.
\begin{proposition}
  \label{theoVB}
  For a standard Bernstein function $v:(0,\infty)\to\mathbb{R}$ and a
  fixed parameter $c>0$, a probability measure $B$ on $(0,1]$ such
  that
  \begin{align}
    \label{BtoV}
    v(t)&=\int_{(0,1]}\frac{(c t+1)^\beta-1}{\beta c}dB(\beta)
  \end{align}
  is unique if it exists.
\end{proposition}
Thus mixtures of power laws can be identified if the mixture exists.

A careful reader notices that the span of power laws in (\ref{BtoV})
is rather narrow. In practical applications, we need an
identifiability result for mixtures with a varying parameter
$c>0$. Solving this problem requires a longer excursion to functional
analysis so we leave the reader in suspension.

\subsection{Five models}
\label{secModels}

Once we have clarified the general theory of vocabulary growth models,
we proceed to concrete examples. We present five Bernstein functions
of a linguistic relevance. Some of these functions can be connected to
well known formulas for Zipf's and Heaps' laws and three of them were
discussed half-formally by \citet{Debowski25}.

For simplicity of formulas, we exhibit either standard functions or
singular functions with $v(0+)=0$ and $v(1)=1$. We couple Bernstein
functions with their corresponding hapax rates and, as proposed by
\citet{Debowski25}, we refer to these pairs as \strong{models}. The
names of models appeal to the shape of function $h(u)$, which is more
distinctive than the plot of function $v(t)$.  We also denote the rank
function
\begin{align}
  r(t|f)
  =
  v(t)-\sum_{k=1}^{f-1} v(t|k),
\end{align}
which captures the rank-frequency plot.

The simplest model, called the maximal model, consists of types that
have the intensity equal $0$.  Consequently, all word tokens are
hapaxes! Obviously, this model is standard and expansive.
\begin{example}
  \label{exMaximal}
  The \strong{maximal model} is given as
  \begin{align}
    \label{MaximalH}
    h(u)&=1,
    \\
    \label{MaximalV}
    v(t)&=t,
    \\
    \label{MaximalVK}
    v(t|k)&=\boole{k=1}t.
  \end{align}
  The maximal model defines an expansive and standard Bernstein
  function.
\end{example}

\begin{proposition}
  \label{theoMaximalF}
  The Bernstein measure for the model in Example \ref{exMaximal} is
  \begin{align}
    \label{MaximalF}
    F(\klam{0})
    &=1,
      \quad f(p)=0.
  \end{align}
\end{proposition}

The next model, called the constant model, assumes that the proportion
of hapaxes does not depend on the text length. This model leads to the
well known Heaps law and is singular.
\begin{example}
  \label{exConstant}
  The \strong{constant model} with a parameter $\beta\in(0,1)$ is
  given as
  \begin{align}
    \label{ConstantH}
    h(u)&=\beta,
    \\
    \label{ConstantV}
    v(t)&=t^\beta,
    \\
    \label{ConstantVK}
    v(t|k)&=-\frac{(0-\beta)(1-\beta)\ldots(k-1-\beta)}{k!}\cdot t^\beta.
  \end{align}
  Relationship (\ref{ConstantV}) is known as Heaps' law or Herdan's
  law in quantitative linguistics \citep{Herdan64,Heaps78}.  The
  corresponding rank function is
  \begin{align}
    \label{ConstantRF}
    r(t|f)
    =
    \frac{(1-\beta)(2-\beta)\ldots(f-1-\beta)}{(f-1)!}\cdot t^\beta
    \xrightarrow[f\to\infty]{} 0,
  \end{align}
  which is close to the Zipf-Mandelbrot law \citep{Mandelbrot54} for
  large frequencies $f$ but differs systematically for small $f$.  The
  constant model defines an expansive and singular Bernstein function.
\end{example}

Let us define the gamma function
\begin{align}
  \Gamma(z)=\int_0^\infty t^{z-1}e^{-t}dt,\quad \Re z>0.
\end{align}
\begin{proposition}
  \label{theoConstantF}
  The Bernstein measure for the model in Example \ref{exConstant} is
  \begin{align}
   \label{ConstantF}
    F(\klam{0})
    &=0,
      \quad f(p)=\frac{\beta}{\Gamma(1-\beta)}\cdot p^{-\beta}.
  \end{align}
\end{proposition}

Despite its simplicity, the constant model is a visibly incorrect
model of real textual data. \citet{Debowski25} observed that the
empirical hapax rate $h(u)$ is a monotonically decreasing function for
novel-sized texts.  This trend is approximately linear in variable
$u$. There are three simple remedies.

The first idea is to apply an affine transformation to the constant
model, which yields the saturation model. This model is standard and
expansive.
\begin{example}
  \label{exSaturation}
  The \strong{saturation model} with a parameter
  $\beta\in(0,1)$ is given as
  \begin{align}
    \label{SaturationH}
    h(u)&=\frac{\beta e^u(e^u+1)^{\beta-1}}{(e^u+1)^{\beta}-1}
          \to
          \begin{cases}
            1, & u\to -\infty,
            \\
            \beta, & u\to \infty,
          \end{cases}
    \\
    \label{SaturationV}
    v(t)&=\beta^{-1}\okra{(t+1)^\beta-1},
    \\
    \label{SaturationVK}
    v(t|k)&=-\frac{(0-\beta)(1-\beta)\ldots(k-1-\beta)}{k!}
            \okra{\frac{t}{t+1}}^k (t+1)^\beta.
  \end{align}
  The saturation model defines an expansive and standard Bernstein
  function.
\end{example}

\begin{proposition}
  \label{theoSaturationF}
  The Bernstein measure for the model in Example \ref{exSaturation} is
  \begin{align}
   \label{SaturationF}
    F(\klam{0})
    &=0,
      \quad f(p)=\frac{1}{\Gamma(1-\beta)}\cdot p^{-\beta}e^{-p}.
  \end{align}
\end{proposition}

Another approach is the parameter-free cancellation model. This model
recovers the exact Zipf law and is expansive but is also singular.
\begin{example}
  \label{exCancellation}
  The \strong{cancellation model} is given as
  \begin{align}
      \label{CancellationH}
    h(u)&=
        \begin{cases}
          \displaystyle
          \frac{1}{u}-\frac{1}{e^{u}-1}, & u\neq 0,
          \\
          \frac{1}{2}, & u=0,
        \end{cases}
    \\
    \label{CancellationV}
    v(t)&=
        \begin{cases}
          \displaystyle
          \frac{t\log t}{t-1}, & t\neq 1,
          \\
          1, & t=1.
        \end{cases}
    \\
    \label{CancellationVK}
    v(t|k)
      &=
        \sum_{j=k}^\infty \binom{j}{k} \frac{t^{k}(1-t)^{j-k}}{j(j+1)},
        \quad t\in[0,2].
  \end{align}
  The poles at $u=0$ in formula (\ref{CancellationH}) cancel, hence
  the name.  The cancellation model was studied by \citet[p.\
  97--101]{Baayen01} and rediscovered by \citet{Davis18}. The
  corresponding rank function is
  \begin{align}
    \label{CancellationRF}
    r(t|f)
    &=
    \begin{cases}
      \displaystyle
      \frac{\log t-\sum_{j=1}^{f-1} \frac{\okra{1-\frac{1}{t}}^{j}}{j}
      }{\okra{1-\frac{1}{t}}^{f}}, & t\neq 1,  
      \\
      \displaystyle
      \frac{1}{f}, & t=1, 
    \end{cases}
  \end{align}
  which yields the exact Zipf law \citep{Zipf35} for $t=1$. For
  $t\ge 1/2$, formula (\ref{CancellationRF}) can be rewritten as
  \begin{align}
    \label{CancellationRFII}
    r(t|f)
    &=\sum_{j=0}^\infty \frac{\okra{1-\frac{1}{t}}^{j}}{j+f}
      \xrightarrow[f\to\infty]{} 0.
  \end{align}
  Formulas (\ref{CancellationVK})--(\ref{CancellationRFII}) were
  derived by \citet[Section 4.2]{Debowski25}. The model defines a
  Bernstein function since $v(t|k)\ge 0$ either by formula
  (\ref{CancellationVK}) for $t\le 1$ or by (\ref{CancellationRFII})
  for $t\ge 1$. This Bernstein function is singular and expansive for
  $t\ge 1$ by (\ref{CancellationRFII}).
\end{example}

The next proposition applies an exponential integral given as
\begin{align}
  E_2(z):=\int_{1}^{\infty} \frac{e^{-zu}du}{u^2},\quad \Re z>0.
\end{align}
\begin{proposition}
  \label{theoCancellationF}
  The Bernstein measure for the model in Example \ref{exCancellation} is
  \begin{align}
   \label{CancellationF}
    F(\klam{0})
    &=0,
      \quad f(p)=e^p E_2(p).
  \end{align}
\end{proposition}

Yet another approximation is given by the logistic model, whose affine
transformation yields a smaller fitting error for novel-sized texts
than that of the cancellation model \citep{Debowski25}. This model was
proposed by \citet{Debowski25} but not fully analyzed theoretically.
\begin{example}
  \label{exLogistic}
  The \strong{logistic model} with a parameter $\gamma>0$ is
  given as
  \begin{align}
    \label{LogisticH}
    h(u)&=\frac{1}{e^{\gamma u}+1},
    \\
    \label{LogisticV}
    v(t)&=\frac{t}{(t^{\gamma}+1)^{1/\gamma}}.
  \end{align}
  In particular, for $\gamma=1$, we obtain the rank function
  \begin{align}
    \label{LogisticRF}
    r(t|f)
    &=\frac{2 t^{f}}{(t+1)^{f}}
      \xrightarrow[f\to\infty]{} 0,
  \end{align}
  Hence $v(t|k)\ge 0$ and the model defines a Bernstein function for
  $\gamma=1$. This Bernstein function is bounded and standard.
\end{example}

Justifying $v(t|k)\ge 0$ in the case of $\gamma\neq 1$ was stated as
an open problem by \citet{Debowski25}.  Now we solve it.
\begin{proposition}
  \label{theoLogistic}
  The model defined in Example \ref{exLogistic} defines a bounded and
  standard Bernstein function for $\gamma\in(0,1]$.
\end{proposition}
Akin to the methods of \citet{Debowski25}, the basic technique in our
proof of Proposition \ref{theoLogistic} is the analysis of subsequent
derivatives of function $v(t)$ by means of a polynomial recursion. As
compared to the previous work, we make some new observations, which
allow us to draw the desired conclusions.

In fact, Proposition \ref{theoLogistic} can be alternatively proved
referring to the stable density $g_\gamma(p;c)$ defined in Example
\ref{exStable}.
\begin{proposition}
  \label{theoLogisticF}
  The Bernstein measure for the model in Example \ref{exLogistic} is
  \begin{align}
   \label{LogisticF}
    F(\klam{0})
    &=0,
      \quad f(p)=\frac{1}{\Gamma(1+1/\gamma)}
      \int_0^\infty c^{1/\gamma}e^{-c}g_\gamma(p;c)dc.
  \end{align}
\end{proposition}

An interested reader is referred to the table of Bernstein functions
and their measures presented by \citet[Chapter
15]{SchillingSongVondracek10}. Combined with examples presented in
this section, they offer some toolbox of potential vocabulary growth
models for a language researcher.

\section{Multinomial sources}
\label{secMultinomial}

The second kind of a source to be analyzed is the multinomial
source. In this case, there also exists a rich theory. The
developments of this section are mostly quoted after the work of
\citet[Appendix A.4]{Debowski25g}.  In particular, the expectations of
our quantities of interest can be computed exactly as the following
proposition assures.
\begin{proposition}
  \label{theoMeanD}
  For the multinomial source, under condition (\ref{Probability}) and
  $n\in\mathbb{N}$, we have
  \begin{align}
    \label{TokensD}
    N_M(n)
    &=n,
    \\
    \label{SmallVD}
    \mean V_M(n)
    &=v_M(n)
    := \sum_{w\in\mathbb{W}} (1-(1-p_w)^{n}),
    \\
    \label{SmallVKD}
    \mean V_M(n|k)
    &=v_M(n|k)
    :=
      \begin{cases}
        \sum_{w\in\mathbb{W}} \binom{n}{k} p_w^{k}(1-p_w)^{n-k},
        & 1\le k\le n,
        \\
        0,
        & k>n.               
      \end{cases}
  \end{align}
\end{proposition}

Now we will see that the expected spectrum can be computed from the
expected number of types by taking differences. In the
following, we write the  \strong{difference operator}
\begin{align}
  \label{DerivativesD}
  \Delta f(n)&:=f(n)-f(n-1).
\end{align}
We also write $\Delta^0 f(n):=f(n)$ and
$\Delta^{k+1} f(n):=\Delta \Delta^k f(n)$.
We assume that this operator acts only on the directly following
functor: $\Delta f(g(n))=f(g(n))-f(g(n)-1)$.
\begin{proposition}
  \label{theoTaylorD}
  For the multinomial source, under $n\in\mathbb{N}$ and $t>0$, we
  have
  \begin{align}
    \label{TaylorD}
    v_M(n|k)
    &=
      \begin{cases}
        (-1)^{k+1}\binom{n}{k}\Delta^k v_M(n),
        & 1\le k\le n,
        \\
        0,
        & k>n.
      \end{cases}
  \end{align}
\end{proposition}

\subsection{Hausdorff sequences}
\label{secHausdorffI}

By an analogy to Section \ref{secBernsteinI}, for arbitrary function
$v_M(n)$, let us treat claim (\ref{TaylorD}) as the definition of
functions $v_M(n|k)$. Consequently, conditions $v_M(n)\ge 0$ and
$v_M(n|k)\ge 0$ imply that function $v_M(n)$ is non-negative and its
subsequent differences have alternating signs. Such a concept is also
well known in mathematics.
\begin{definition}[\mbox{\citealp{Akhiezer21}}]
  A sequence $v:\mathbb{N}\cup\klam{0}\to\mathbb{R}$ is called a
  \strong{Hausdorff sequence} if $v(n)\ge 0$ and
  $(-1)^{k+1}\Delta ^k v(n+k)\ge 0$ for all $k\in\mathbb{N}$ and
  $n\in\mathbb{N}\cup\klam{0}$.  A sequence
  $u:\mathbb{N}\cup\klam{0}\to\mathbb{R}$ is called \strong{completely
    alternating} if $u(n)=\Delta v(n+1)$ for a certain Hausdorff
  sequence $v:\mathbb{N}\cup\klam{0}\to\mathbb{R}$.
\end{definition}
We call a Hausdorff sequence $v(n)$ \strong{standard} if $v(0)=0$ and
$\Delta v(1)=1$. The name \emph{Hausdorff sequence} is non-standard
itself and coined by an analogy to the standard term \emph{Bernstein
  function}, defined in Section \ref{secBernsteinI}.

Let us present some simple examples of Hausdorff sequences.
\begin{example}
  Obviously, function $v_M(n)$ defined in (\ref{SmallVD}) is a
  Hausdorff sequence.  Here are some more examples of Hausdorff
  functions: $n\mapsto n$, $n\mapsto \sum_{m=1}^n m^{-1}$, and
  $n\mapsto p^{-1}\okra{1-(1-p)^n}$ for and $p>0$.  All these examples
  are standard.
\end{example}

\subsection{Hausdorff theorem}
\label{secHausdorffII}

It is true that any Hausdorff sequence is a convex combination of
sequences $n\mapsto p^{-1}\okra{1-(1-p)^n}$ for varying $p>0$. A
particular case of a discrete combination has been already encountered
in formula (\ref{SmallVD}). The more general result is known as the
Hausdorff moment theorem \citep{Hausdorff23}.
\begin{proposition}[\mbox{\citealp{Akhiezer21}}]
  \label{theoVFD}
  A sequence $v:\mathbb{N}\cup\klam{0}\to\mathbb{R}$ is a Hausdorff
  sequence if and only if there exists a unique non-negative measure
  $F$ on $[0,1]$ such that and
  \begin{align}
    \label{FtoVD}
    v(n)&=v(0)
          +
          nF(\klam{0})
          +
          \int_{(0,1]}\frac{1-(1-p)^n}{p}dF(p),
    \\
    \label{FtoVPrimeD}
    \Delta v(n+1)
        &=F(\klam{0})
          +
          \int_{(0,1]}(1-p)^n dF(p).
  \end{align}
\end{proposition}
Additionally, we see that a Hausdorff sequence
$v:\mathbb{N}\cup\klam{0}\to\mathbb{R}$ is standard if and only if
$v(0)=0$ and the measure $F$ is a \strong{probability measure},
namely,
\begin{align}
  F(\klam{0})+\int_{(0,1]} dF(p)=1.
\end{align}
We call measure $F$ the \strong{Hausdorff measure}.
\begin{example}
  The Hausdorff sequence $v(n)=v_M(n)$ defined in (\ref{SmallVD}) has a
  Hausdorff measure that is supported on discrete points
  \begin{align}
    \label{Points}
    F(A)=\sum_{w\in\mathbb{W}} p_w\boole{p_w\in A}, \quad p_w\le 1.
  \end{align}
\end{example}

Many useful Hausdorff sequences enjoy a continuous Hausdorff measure that
has a density $f:(0,\infty)\to\mathbb{R}$ with respect to the Lebesgue
measure.  We call $f$ the \strong{Hausdorff density}.  Then we may
substitute
\begin{align}
  dF(p)=f(p)dp.
\end{align}
\begin{example}
  Let the Hausdorff probability density be $f(p)=1$.  We obtain a
  standard Hausdorff sequence $v(n)$ that has the first difference
  \begin{align}
    \Delta v(n+1)
    &=
      \int_{[0,1]}(1-p)^ndp=\frac{1}{n+1}.
  \end{align}
  Hence $v(n)=\sum_{m=1}^n\Delta v(m)$ grows asymptotically as
  $\log n$.
\end{example}

To produce a vocabulary growth that is closer to a power law, we may
consider the following class of Hausdorff measures.
\begin{example}
  Let the Hausdorff probability density be $f(p)=(1-\beta)p^{-\beta}$,
  where $\beta\in(0,1)$.  We obtain a standard Hausdorff sequence
  $v(n)$ that has the first difference
  \begin{align}
    \Delta v(n+1)
    &=
      (1-\beta)\int_{[0,1]}(1-p)^n p^{-\beta}dp
      =
    \frac{\Gamma(n+1)\Gamma(2-\beta)}{\Gamma(n+2-\beta)}
    \nonumber\\
    &=
      \frac{n!}{(2-\beta)(3-\beta)\ldots(n+1-\beta)}\sim n^{1-\beta}.
  \end{align}
  Hence $v(n)=\sum_{m=1}^n\Delta v(m)$ grows asymptotically as
  $n^\beta$.
\end{example}

\subsection{Newton formula}
\label{secNewton}

For an arbitrary sequence, let us define spectrum elements as in
(\ref{TaylorD}).
\begin{definition}
  For a sequence $v:\mathbb{N}\cup\klam{0}\to\mathbb{R}$, we
  define the \strong{spectrum elements}
  \begin{align}
  \label{TaylorIID}
  v(n|k)&:=(-1)^{k+1}\binom{n}{k} \Delta^k v(n), &  1\le k\le n.
  \end{align}
\end{definition}

We notice that the familiar Newton formula
$(a+b)^m=\sum_{k=0}^m \binom{m}{k} a^kb^{m-k}$ can be rewritten in the
following form, which resembles the Taylor expansion discussed in
Section \ref{secTaylor}.
\begin{proposition}
  \label{theoNewton}
  For a sequence $v:\mathbb{N}\cup\klam{0}\to\mathbb{R}$ and
  $0\le m\le n$, we have
  \begin{align}
    \label{Newton}
    v(n-m)
    &=
      \sum_{k=0}^m (-1)^k
      \binom{m}{k} \Delta^k v(n).
  \end{align}
\end{proposition}

By formula (\ref{Newton}), consistency conditions analogous to
(\ref{Consistent}) can be asserted easily for general sequences.
\begin{proposition}
  \label{theoTypeConsistentD}
  A sequence $v:\mathbb{N}\cup\klam{0}\to\mathbb{R}$
  satisfies condition $\sum_{k=1}^n v(n|k)=v(n)$ if and only if
  $v(0)=0$.
\end{proposition}
\begin{proposition}
  \label{theoTokenConsistentD}
  A sequence $v:\mathbb{N}\cup\klam{0}\to\mathbb{R}$ satisfies
  condition $\sum_{k=1}^n kv(n|k)=n$ if and only if
  $\Delta v(1)=1$.
\end{proposition}
Thus only standard Hausdorff sequences satisfy both consistency
conditions and can be applied as models of word frequency
distributions.

\subsection{Hapax rate}
\label{secHapaxD}

Another perspective onto Hausdorff sequences is given by considering
the proportion of types that occur exactly once.
\begin{definition}
  For a Hausdorff sequence $v:\mathbb{N}\cup\klam{0}\to\mathbb{R}$, we
  define the \strong{hapax rate}
  $h:\mathbb{N}\setminus\klam{1}\to[0,1]$ where
  \begin{align}
    \label{RateD}
    h(n)&:=\frac{v(n|1)}{v(n)-v(0)}.
  \end{align}
\end{definition}

According to the next statement, there is a one-to-one correspondence
between hapax rates and Hausdorff sequences if we fix the value of the
latter at arguments $n\in\klam{0,1}$.
\begin{proposition}
  \label{theoVHD}
  For a sequence $h:\mathbb{N}\setminus\klam{1}\to[0,1]$ and a
  sequence $v:\mathbb{N}\cup\klam{0}\to\mathbb{R}$, condition
  (\ref{RateD}) holds if and only if
  \begin{align}
    \label{HtoVD}
    v(n)-v(0)&=(v(1)-v(0))\prod_{m=2}^{n}\okra{1-\frac{h(m)}{m}}^{-1}.
  \end{align}
\end{proposition}
Not every sequence $h:\mathbb{N}\setminus\klam{1}\to[0,1]$ is a hapax
rate of a certain Hausdorff sequence but each valid hapax rate $h(n)$
defines exactly one standard Hausdorff sequence
\begin{align}
  v(n)=\prod_{m=2}^{n}\okra{1-\frac{h(m)}{m}}^{-1}.
\end{align}

\section{Invariant sources}
\label{secInvariant}

In the general case of stationary processes, the expected number of
types need not be completely alternating but in some properties, it
resembles the Shannon entropy.  For a stationary process
$(X_i)_{i\in\mathbb{Z}}$, we denote the block Shannon entropy and the
expected number of types as
\begin{align*}
  H(n)&:=\mean\okra{-\log P(X_1,X_2,\ldots,X_n)},
  \\
  v(n)&:=\mean\card\klam{X_1,X_2,\ldots,X_n}.
\end{align*}

The relevant results from the work of \citet{Debowski25g} are as
follows.  Function $f(n)$ is called positive, growing, and concave if
$f(n)\ge 0$, $\Delta f(n)\ge 0$, and $\Delta^2 f(n)\le 0$,
respectively.
\begin{proposition}
  \label{theoBlockEntropy}
  Let a stationary process $(X_i)_{i\in\mathbb{Z}}$.  For
  $n\in\mathbb{N}$, we claim that:
  \begin{enumerate}   
  \item Function $n\mapsto H(n)$ is positive, growing, and concave and
    $H(0)=0$.
  \item Function $n\mapsto H(n)/n$ is decreasing.
  \item There exists limit $\lim_{n\to\infty}H(n)/n\ge 0$.
  \end{enumerate}  
\end{proposition}
\begin{proposition}
  \label{theoBlockTypes}
  Let a stationary process $(X_i)_{i\in\mathbb{Z}}$.  For
  $n\in\mathbb{N}$, we claim that:
  \begin{enumerate}   
  \item Function $n\mapsto v(n)$ is positive, growing, and concave and
    $v(0)=0$.
  \item Function $n\mapsto v(n)/n$ is decreasing and $v(n)/n\le 1$.
  \item There exists limit $\lim_{n\to\infty} v(n)/n=0$.
  \end{enumerate}  
\end{proposition}
For some dependent processes (such as periodic processes), the second
differences $\Delta^2 H(n)$ and $\Delta^2 v(n)$ may oscillate as a
function of $n$.

In the stationary case, there is a general upper bound for the
expected number of hapaxes and dis legomena (twice-comers) in terms of
the difference of the total number of types.
\begin{proposition}
  \label{theoHapaxes}
  For a stationary process $(X_i)_{i\in\mathbb{Z}}$, we have
  \begin{align}
    \label{Hapaxes}
    \frac{v(n|1)}{n}\le
    \Delta v\okra{\ceil{\frac{n}{2}}}.
  \end{align}
\end{proposition}
\begin{proposition}
  \label{theoDis}
  For a stationary process $(X_i)_{i\in\mathbb{Z}}$, we have
  \begin{align}
    \label{Dis}
    \frac{v(n|2)}{n}\le
    -\sum_{m=1}^{n-1} m\Delta^2 v(n-m+1).
  \end{align}
\end{proposition}
Formulas (\ref{Hapaxes})--(\ref{Dis}) are analogues of
formula (\ref{TaylorD}).

\section{Weibull sources}
\label{secWeibull}

For Weibull sources, exact results are somewhat limited. This section
is mostly an invitation to future research.  We note that if the
generating process exhibits burstiness then hapaxes tend to occur near
the edge of the text.
\begin{proposition}
  \label{theoWeibull}
  For the Weibull source with $\alpha_w\le 1$, under $t>0$, we have
  \begin{align}
    \label{SmallVW}
    \mean V_W(t)
    &=\sum_{w\in\mathbb{W}}
      \okra{1-e^{-(t/\lambda_w)^{\alpha_w}}},
    \\
    \mean V_W(t|k)
    &\le\sum_{w\in\mathbb{W}}
      \frac{\Gamma(\alpha_w+1)^k}{\Gamma(k\alpha_w+1)}
      \okra{\frac{t}{\lambda_w}}^{k\alpha_w}e^{-(t/\lambda_w)^{\alpha_w}},
    \\
    (-1)^{k+1}t^k\nabla^k
    \mean V_W(t)
    &\ge\sum_{w\in\mathbb{W}}
    \alpha_w^k \okra{\frac{t}{\lambda_w}}^{k\alpha_w}e^{-(t/\lambda_w)^{\alpha_w}}.
  \end{align}
\end{proposition}
We notice that function $\mean V_W(t)$ calculated in (\ref{SmallVW})
is a Bernstein function by a reasoning applying Example
\ref{exStable}.
 
Function $\alpha\mapsto \frac{\Gamma(\alpha)^k}{\Gamma(k\alpha+1)}$ is
monotonically decreasing for $0\le\alpha\le 1$.  Hence, by Proposition
\ref{theoWeibull}, if we have a uniform lower bound on burstiness
$\alpha<\alpha_w\le 1$ for $w\in\mathbb{W}$ then we obtain an analogue
of formula (\ref{Taylor}), namely,
\begin{align}
  \mean V_W(t|k)
  &\le
    \frac{\Gamma(\alpha)^k}{\Gamma(k\alpha+1)}
    (-1)^{k+1} t^k\nabla^k
    \mean V_W(t).
\end{align}
We recall that usually $0\le\alpha_w\le 1$ for words of natural
language.

\section{Conclusion}
\label{secConclusion}

In this paper, we have offered a glimpse onto a mathematical theory of
Bernstein functions and Hausdorff sequences from a quantitative
linguistic perspective. We aimed to sketch the theoretical breadth of
the simple topic that starts with counting distinct words in a
text. Our guiding idea has been that marginal lexical statistics,
while heuristically easy to describe at the level of Zipf's and Heaps'
laws, admit a richer mathematical structure.  By connecting linguistic
questions to a well developed branch of analysis, we hope to
contribute to a more rigorous foundation for quantitative
linguistics. We also hope that this work will stimulate further
dialogue between mathematicians and linguists.



\section*{Funding}

This work received no funding from funding agencies.

\section*{Disclosure statement}

There are no competing interests to declare.

\section*{AI usage}

The general idea of this paper originated from a dialogue of the
author with ChatGPT prompted about completely monotone functions and
completely alternating sequences. The chatbot was used as an
intelligent search engine for more established mathematical
notions. In particular, it drew the author's attention to the concept
of Bernstein functions, the Hausdorff moment theorem, and stable
densities, providing initial references to the existing
literature. Except for some paragraphs drafted by the chatbot in the
Introduction, the paper organization, the selection of described
facts, and the previously unpublished proofs are due to the author.



\appendix

\section{Proofs for Poisson sources}
\label{secPoissonA}

In this appendix, we demonstrate Propositions
\ref{theoMean}--\ref{theoLogisticF} with some exceptions.  The proofs
of Propositions \ref{theoMean}--\ref{theoTaylor} are elementary. The
less obvious parts of them are hinted, the more evident ones are left
as exercises. For proofs of Propositions
\ref{theoVF}--\ref{theoAnalyticV}, we refer to
\citet{SchillingSongVondracek10}.  Propositions
\ref{theoPowerSeries}--\ref{theoTokenConsistent} follow easily from
Proposition \ref{theoAnalyticV} but we could not locate a suitable
reference, so we demonstrate them explicitly.  For proofs of
Propositions \ref{theoVH}--\ref{theoH} we refer to \citet{Debowski25}
and \citet{Baayen01}.  The proofs of Propositions \ref{theoHStandard}
and \ref{theoBoundedExpansive} seem novel.  For the proof of
Proposition \ref{theoChoquet}, we refer to \citet{Phelps01}.  The
elementary proof of Proposition \ref{theoVB} is given.  The proof of
Propositions \ref{theoMaximalF}--\ref{theoLogisticF} follow by
arguments known in the literature of Bernstein functions
\citep{SchillingSongVondracek10} but we reproduce them for
completeness.

\subsection{Proof of Proposition \ref{theoMean}}

Assume distribution (\ref{Poisson}).  Equations (\ref{Tokens}) and
(\ref{VarTokens}) follow easily since the sum of Poisson distributed
random variables is Poisson distributed with the intensity equal to
the sum of contributing intensities.

Next, we write $N^P_w:=N^P_w((0,t])$.  We have
$P(N^P_w>0)=1-e^{-p_w t}$ and
\begin{align}
  P(N^P_w>0)\ge P(N^P_w>0)(1-P(N^P_w>0))= 1-e^{-2p_w t}-1+e^{-p_w t}.
\end{align}
Thus, by independence, we obtain
\begin{align}
  \mean V_P(t)&=\sum_w \mean\boole{N^P_w>0}
            =\sum_w P(N^P_w>0)=v_P(t),
  \\
  \var V_P(t)&=\sum_w \var\boole{N^P_w>0}
               \nonumber\\
              &=\sum_w P(N^P_w>0)(1-P(N^P_w>0))
                \nonumber\\
              &=v_P(2t)-v_P(t)\le v_P(t).
\end{align}

Similarly, we obtain
\begin{align}
  P(N^P_w=k)(1-P(N^P_w=k))= \frac{(p_w t)^k}{k!} e^{-p_w t}
  -\frac{1}{2^{2k}}\binom{2k}{k}\frac{(2p_w t)^{2k}}{(2k)!} e^{-2p_w t}.
\end{align}
Hence, analogously, we derive
\begin{align}
  \mean V_P(t|k)&=\sum_w \mean\boole{N^P_w=k}
            =\sum_w P(N^P_w=k)=v_P(t|k),
  \\
  \var V_P(t|k)&=\sum_w \var\boole{N^P_w=k}
             \nonumber\\
            &=\sum_w P(N^P_w=k)(1-P(N^P_w=k))
              \nonumber\\
            &=v_P(t|k)-\frac{1}{2^{2k}}\binom{2k}{k}v_P(2t|2k)\le v_P(t|k).
\end{align}

\subsection{Proof of Proposition \ref{theoTaylor}}

The claims follow by the simple observation $\nabla^k e^{-p t}=p^k e^{-p t}$.

\subsection{Proof of Proposition \ref{theoPowerSeries}}

Let a Bernstein function $v(t)$.  By property (\ref{Taylor}), the
right-hand side of formula (\ref{PowerSeries}) can be written as the
Taylor series
\begin{align}
  v_t(z)=v(t)+\sum_{k=1}^\infty\frac{(-z t)^k}{k!}\nabla^k v(t).
\end{align}
Since the analytic extension $v(z)$ of $v(t)$ exists on the
half-plane $\Re z>0$, we have $v_t(z)=v(t-z t)$ in the open disk
$\abs{z}<1$ \citep[Proposition 5.18]{Lvovski20}.  Consider now a real
number $s\in[0,1]$. Since the spectrum elements are non-negative, by
the monotone convergence theorem \citep[Theorem 16.2]{Billingsley79}
and continuity of function $v(t-s t)$ at point $s=1$, we obtain
\begin{align}
  v(0+)
  =\lim_{s\to 1}v(t-s t)
  =\lim_{s\to 1}v_t(s)
  =v(t)-\sum_{k=1}^\infty v(t|k)
\end{align}
Because $v(0+)$ is finite, we derive
$\sum_{k=1}^\infty v(t|k)=v(t)<\infty$.  Hence, since the spectrum
elements $v(t|k)$ are non-negative, series $v_t(z)$ converges
absolutely for all $\abs{z}\le 1$. Consequently, function $v_t(z)$ is
continuous for all $\abs{z}\le 1$. By continuity of $v(z)$ for
$\Re z>0$, we infer $v_t(z)=v(t-z t)$ also on the unit circle
$\abs{z}=1$.

Subsequently, since the analytic extension $v(z)$ of $v(t)$ exists on
the half-plane $\Re z>0$, we can differentiate series $v_t(z)$ term by
term with respect to variable $z$ for $\abs{z}<1$. In this way, we
obtain
\begin{align}
  t\nabla v(t-z t)=
  -\nabla v_t(z)=\sum_{k=1}^\infty kz^{k-1} v(t|k)
\end{align}
for $\abs{z}<1$.  Taking the limit $z\to 1$, we obtain
$\sum_{k=1}^\infty k v(t|k)=t\nabla v(0+)$, regardless whether $\nabla v(0+)$ is
finite, by the monotone convergence theorem \citep[Theorem
16.2]{Billingsley79} since the spectrum elements $v(t|k)$ are
non-negative.

\subsection{Proof of Proposition \ref{theoTypeConsistent}}

The claim follows by formula (\ref{PowerSeries}) for $z=1$.

\subsection{Proof of Proposition \ref{theoTokenConsistent}}

The claim follows by formula (\ref{PowerSeriesPrime}) for $z=1$.

\subsection{Proof of Proposition \ref{theoHStandard}}

Differentiating formula (\ref{HtoV}), we obtain
\begin{align}
  \nabla v(t)
  &=v(1) \exp\okra{\int_0^{\log t} h(u)du}
    h(\log t) \frac{d\log t}{dt}.
    \nonumber\\
  &=v(1) \exp\okra{\int_0^{\log t} h(u)du}
    h(\log t) \frac{1}{t}
    \nonumber\\
  &=v(1) \exp\okra{-\int_0^{\log t}(1-h(u))du}
    h(\log t).
\end{align}
Hence $\nabla v(0+)<\infty$ if and only if (\ref{RetroDampedH}) is true
because $\lim_{u\to -\infty}h(u)=1$ holds whenever
(\ref{RetroDampedH}) holds. Hence, $\nabla v(0+)=1$ if and only if we have
(\ref{V1}).

\subsection{Proof of Proposition \ref{theoBoundedExpansive}}

By contradiction, let a Bernstein function that is both bounded and
expansive at all points $t> 0$. We can evaluate the derivative
\begin{align}
  \nabla v(t|k)
  &=
  -\frac{(-1)^kn^{k-1}}{(k-1)!}\nabla^k v(t)
    -\frac{(-1)^kn^k}{k!}\nabla^{k+1} v(t)
    \nonumber\\
  &=
    \frac{kv(t|k)-(k+1)v(t|k+1)}{t}.
    \label{DerivativeVK}
\end{align}
Hence, since the spectrum elements $v(t|k)$ are non-negative and the
Bernstein function is expansive at all $t\ge 0$, we have
\begin{align}
  \label{BoundedVK}
  \frac{v(t|k+1)}{v(t|k)}&\le\frac{k}{k+1}
\end{align}
for all $k\ge 1$. In view of formula (\ref{PowerSeries}) for $z=-1$
and the upper bound (\ref{BoundedVK}), we may write
\begin{align}
  v(2 t)-v(t)
  &=
    \sum_{k=1}^\infty (-1)^{k+1} v(t|k)
    \nonumber\\
  &=
    \sum_{k=1}^\infty (v(t|2k-1)-v(t|2k))
    \nonumber\\
  &\ge
    \sum_{k=1}^\infty v(t|2k-1)\okra{1-\frac{2k-1}{2k}}
    =
    \sum_{k=1}^\infty \frac{v(t|2k-1)}{2k}.
  \label{Alternating}
\end{align}
Because the Bernstein function is bounded, we have
$\lim_{t\to\infty} v(t)<\infty$.  Hence
\begin{align}
  \lim_{t\to\infty} (v(2 t)-v(t))=0.
\end{align}
Applying this to (\ref{Alternating}) and invoking non-negative
spectrum elements $v(t|k)$ and the upper bound (\ref{BoundedVK})
again, we obtain
\begin{align}
  \lim_{t\to\infty} v(t|k)=0
\end{align}
for all $k\ge 1$. On the other hand, because spectrum elements
$v(t|k)$ are non-negative and they sum up to $v(t)$, we also derive
$v(1|k)>0$ for some $k\ge 1$. This means that $\nabla v(t|k)<0$ for some
$t\ge 1$ and $k\ge 1$. But this contradicts our assumption that the
the Bernstein function is expansive at all $t>0$. Hence a Bernstein
function cannot be both closed and expansive at all points.

\subsection{Proof of Proposition \ref{theoVB}}

Differentiating (\ref{BtoV}), we obtain
\begin{align}
  \nabla v(t)=\int_{(0,1]}(c t+1)^{-(1-\beta)}dB(\beta).
\end{align}
Now we will use variable $u=\log(c t+1)$, function
$g(u)=\nabla v(c^{-1}(e^{-u}-1))$, and measure $C(A)=B(\klam{1-p:p\in
  A})$. This yields
\begin{align}
  g(u)
  &=\nabla v(c^{-1}(e^{-u}-1))=\int_{(0,1]}e^{-(1-\beta)u}dB(\beta)
   =\int_{[0,1)}e^{-up}dC(p).
\end{align}
By Proposition \ref{theoVF}, measure $C$ is unique for function
$g(u)$---if it exists. Hence follows uniqueness of measure $B$.

\subsection{Proof of Proposition \ref{theoMaximalF}}

The claim is obvious.

\subsection{Proof of Proposition \ref{theoConstantF}}

We may verify
\begin{align}
  \int_0^\infty e^{-p t} f(p)dp
  =
   \frac{\beta}{\Gamma(1-\beta)}\int_0^\infty p^{-\beta} e^{-p t} dp
  =
  \beta t^{\beta-1}
  =
  \nabla v(t).
\end{align}

\subsection{Proof of Proposition \ref{theoSaturationF}}

We may verify
\begin{align}
  \int_0^\infty e^{-p t} f(p)dp
  =
   \frac{1}{\Gamma(1-\beta)}\int_0^\infty p^{-\beta} e^{-p}e^{-p t} dp
  =
  (t+1)^{\beta-1}
  =
  \nabla v(t).
\end{align}

\subsection{Proof of Proposition \ref{theoCancellationF}}

Let us compute the derivative of $v(t)$ and use its integral
representation
\begin{align}
  \nabla v(t)=\frac{t-1-\log t}{(t-1)^2}=
  \int_{0}^{1}\frac{s}{1+(t-1)s} ds=
  \int_{1}^{\infty}\frac{1}{u-1+t} \frac{du}{u^2}
  .
\end{align}
The integrated expression is the following Laplace transform
\begin{align}
  \frac{1}{u-1+t}
  =\int_{0}^{\infty} e^{-p t}e^{p}e^{-up} dp.
\end{align}
Hence, interchanging the order of integrals by the Fubini theorem
\citep[Theorem 18.3]{Billingsley79}, we obtain
\begin{align}
  \nabla v(t)
  =
  \int_{1}^{\infty}\okra{\int_{0}^{\infty} e^{-p t}e^{p}e^{-up}dp} \frac{du}{u^2}
  =
  \int_{0}^{\infty} e^{-p t}e^{p}\okra{\int_{1}^{\infty}\frac{e^{-up}du}{u^2}} dp. 
\end{align}
Thus, by the uniqueness of the Bernstein measure, we obtain its density
\begin{align}
  f(p)
  =e^p\int_{1}^{\infty} \frac{e^{-up}du}{u^2}=e^p E_2(p).
\end{align}

\subsection{Proof of Proposition \ref{theoLogistic}}

Throughout this section, we fix a parameter $\gamma\in(0,1]$.  All
other conditions being easily met, it remains to prove that the relative
spectrum elements $h(u|k)$ for the logistic model (\ref{LogisticH})
are non-negative.

To begin we observe that
\begin{align}
  \frac{dh(u)}{du}=-\gamma h(u)(1-h(u)).
\end{align}
Hence the relative spectrum elements can be conveniently written as
\begin{align}
  \label{UtoHK}
  h(u|k)&=U_k(1-h(u)),
          \quad k\ge 0,
\end{align}
and the recursion (\ref{HK1})--(\ref{HK2}) can rewritten as
\begin{align}
  U_0(x)
     &=-1,
  \\
  \label{RecursionU}
  U_{k+1}(x)
     &=\frac{k-1+x}{k+1}U_k(x)-\frac{\gamma x(1-x)}{k+1}\nabla U_k(x),
\end{align}
with variable $x=1-h(u)\in[0,1]$ and $\nabla U_k(x)$ being the derivative of
$U_k(x)$ with respect to $x$.

Functions $U_k(x)$ are polynomials of variable $x\in[0,1]$.
Subsequently, we observe that
\begin{align}       
  U_1(x)&=1-x\ge 0,
  \\
  U_2(x)&=\frac{x}{2}(1-x)-\frac{\gamma x(1-x)}{2}(-1) 
       =\frac{1+\gamma}{2} x(1-x)\ge 0.
\end{align}
Factor $x(1-x)$ persists in further iterations.  We can simplify
recursion (\ref{RecursionU}) by considering auxiliary polynomials
$W_k(x)$ defined via
\begin{align}
  \label{WtoU}
  U_k(x)&=\frac{1+\gamma}{2} x(1-x) W_k(x),
          \quad k\ge 2.
\end{align}
Thus we have
\begin{align}
  W_2(x)&=1,
  \\
  \label{RecursionW}
  W_{k+1}(x)&=\frac{k-1+x}{k+1}W_k(x)
              -\frac{d}{dx}\okra{\frac{\gamma x(1-x)}{k+1}W_k(x)}.
\end{align}

Let us explore the properties of recursion (\ref{RecursionW}).  We
will show that polynomials $W_k(x)$ are non-negative, growing, and
convex for $x\ge 0$.  In particular, writing
\begin{align}
  W_k(x)=\sum_{j=0}^{k-2}W_{kj} x^j,
\end{align}
we obtain
\begin{align}
  \label{RecursionWKJ}
  \sum_{j=0}^{k-1}W_{k+1,j} x^j
        &=\sum_{j=0}^{k-2}\frac{W_{kj}x^j}{k+1}
          \okra{k-1+x
          -\gamma(j+1)+\gamma(j+2)x},
\end{align}
Thus, coefficients $W_{kj}$ satisfy the recursion
\begin{align}
  \label{RecursionWKJII}
  W_{k+1,j}&= \frac{1+\gamma(j+1)}{k+1}W_{k,j-1}+
             \frac{k-1-\gamma(j+1)}{k+1}W_{kj}.
\end{align}
Since $W_{kj}=0$ for $j+1>k-1$ and $\gamma\in(0,1]$ then by induction
on $k$, all coefficients $W_{kj}$ are non-negative.  Thus polynomials
$W_k(x)$ are non-negative, growing, and convex for $x\ge 0$.  In view
of formulas (\ref{UtoHK}) and (\ref{WtoU}), this means in particular
that the logistic model has non-negative spectrum elements.

\subsection{Proof of Proposition \ref{theoLogisticF}}

Let us compute the derivative of $v(t)$ and use its integral
representation
\begin{align}
  \nabla v(t)=\frac{1}{(t^\gamma+1)^{1+\frac{1}{\gamma}}}=
  \frac{1}{\Gamma(1+1/\gamma)}
  \int_0^\infty c^{1/\gamma} e^{-c}e^{-c t^\gamma}dc
\end{align}
since for $\beta>0$ we have
\begin{align}
  \frac{1}{(1+u)^{\beta}}
  =\frac{1}{\Gamma(\beta)}\int_0^\infty c^{\beta-1}e^{-c}e^{-cu}dc.
\end{align}
Using expression (\ref{Stable}), we may further write
\begin{align}
  \nabla v(t)=
  \frac{1}{\Gamma(1+1/\gamma)}
  \int_0^\infty c^{1/\gamma} e^{-c}\okra{\int_0^\infty e^{-p t}g_\gamma(p;c)dp}dc.
\end{align}
Thus, interchanging the order of integrals by the Fubini theorem
\citep[Theorem 18.3]{Billingsley79}, we obtain
\begin{align}
  \nabla v(t)
  =
  \int_0^\infty e^{-p t}
  \frac{1}{\Gamma(1+1/\gamma)}\okra{\int_0^\infty c^{1/\gamma}e^{-c}g_\gamma(p;c)dc}dp
  . 
\end{align}
Hence, by the uniqueness of the Bernstein measure, we obtain its density
\begin{align}
  f(p)
  =\frac{1}{\Gamma(1+1/\gamma)}\int_0^\infty c^{1/\gamma}e^{-c}g_\gamma(p;c)dc.
\end{align}


\section{Proofs for multinomial sources}
\label{secMultinomialA}

In this appendix, we demonstrate Propositions
\ref{theoMeanD}--\ref{theoVHD} with some exceptions.  For the proof of
Proposition \ref{theoVFD}, we refer to \citet{Akhiezer21}. Other
results have been mostly derived in the work of \citet[Appendix
A.4]{Debowski25g} but we reproduce their proofs here.

\subsection{Proof of Proposition \ref{theoMeanD}}

Assume distribution (\ref{Binomial}).  Equation (\ref{TokensD})
follows by (\ref{Counting}).  Next, we write $N^M_w:=N^M_w((0,n])$. We
obtain
\begin{align}
  \mean V_M(n)&=\sum_w \mean\boole{N^M_w>0}
            =\sum_w P(N^M_w>0)=v_M(n).
  \\
  \mean V_M(n|k)&=\sum_w \mean\boole{N^M_w=k}
            =\sum_w P(N^M_w=k)=v_M(n|k).
\end{align}.

\subsection{Proof of Proposition \ref{theoTaylorD}}

The claims follow by the simple observation
$\Delta^k(1-p)^n=(-p)^k(1-p)^{n-k}$.

\subsection{Proof of Proposition \ref{theoNewton}}

Observe that formula (\ref{Newton}) follows easily if we denote
$Tf(n)=f(n-1)$.  Then
\begin{align}
  \sum_{k=0}^m (-1)^k
  \binom{m}{k} \Delta^k v(n)
  &=
    \sum_{k=0}^m (-1)^k
    \binom{m}{k} (1-T)^k v(n)
    \nonumber\\
  &=
    \sum_{k=0}^m 
    \binom{m}{k} (T-1)^k 1^{m-k} v(n)
    \nonumber\\
  &=
    (T-1+1)^m v(n)
    =
    T^n v(n)
    =
    v(n-m)
    .
\end{align}

\subsection{Proof of Proposition \ref{theoTypeConsistentD}}

The claim follows by formula (\ref{Newton}) if we put $n=m$ since
\begin{align}
  v(n)-\sum_{k=1}^n v(n|k)
  =
  \sum_{k=0}^n (-1)^k
  \binom{n}{k} \Delta^k v(n)
  =
  v(n-n)
  =
  v(0)
  .
\end{align}

\subsection{Proof of Proposition \ref{theoTokenConsistentD}}

The claim follows by formula (\ref{Newton}) if we put $n=m$ and use
identity $k\binom{n}{k}=n\binom{n-1}{k-1}$ since
\begin{align}
  \sum_{k=1}^n kv(n|k)
  &=
  \sum_{k=1}^n (-1)^{k+1}
  k\binom{n}{k} \Delta^k v(n)
  =
  n\sum_{k=1}^n (-1)^{k+1}
    \binom{n-1}{k-1} \Delta^k v(n)
    \nonumber\\
  &=
  n\sum_{j=0}^{n-1} (-1)^j\binom{n-1}{j} \Delta^{j+1} v(n)
  =
  n\Delta v(n-n+1)
  =
  n\Delta v(1)
  .
\end{align}

\subsection{Proof of Proposition \ref{theoVHD}}

Assume (\ref{RateD}). Then
\begin{align}
  \prod_{m=2}^{n}\okra{1-\frac{h(m)}{m}}^{-1}
  &=
    \prod_{m=2}^{n}\okra{1-\frac{\Delta v(m)}{v(m)-v(0)}}^{-1}
    \nonumber\\
  &=
    \prod_{m=2}^{n}\frac{v(m)-v(0)}{v(m-1)-v(0)}
    =
    \frac{v(n)-v(0)}{v(1)-v(0)}
    .
\end{align}
On the other hand, if we assume (\ref{HtoVD}), we obtain
\begin{align}
  \frac{v(n|1)}{v(n)-v(0)}
  &=
    \frac{n\Delta v(n)}{v(n)-v(0)}
    =
    \frac{n\Delta\prod_{m=2}^{n}\okra{1-\frac{h(m)}{m}}^{-1}}{\prod_{m=2}^{n}\okra{1-\frac{h(m)}{m}}^{-1}}
    \nonumber\\
  &=n\okra{1-1+\frac{h(n)}{n}}=h(n)
    .
\end{align}


\section{Proofs for invariant sources}
\label{secInvariantA}

The proofs of Propositions \ref{theoBlockEntropy}--\ref{theoHapaxes}
can be found in the work of \citet[Appendix
A.3]{Debowski25g}. Proposition \ref{theoDis} seems to be new.

\subsection{Proof of Proposition \ref{theoDis}}

We have
\begin{align}
  v(n|2)
  &=
    \sum_{1\le i<j\le n}
    P(X_i=X_j\not\in X_1^{i-1}
    \cup X_{i+1}^{j-1}\cup X_{j+1}^n)
    .
\end{align}
Hence, we may bound
\begin{align}
  v(n|2)
  &\le
    \sum_{1\le i<j\le n}
    \min\klam{P\okra{X_i\not\in\mathcal{X}_1^{i-1}},
    P\okra{X_i=X_j\not\in\mathcal{X}_{i+1}^{j-1}},
    P\okra{X_j\not\in\mathcal{X}_{j+1}^n}}
    \nonumber\\
  &=
    \sum_{1\le i<j\le n}
    \min\klam{\Delta v(i),
    -\Delta^2 v(j-i+1),
    \Delta v(n-j+1)}
    \nonumber\\
  &\le
    -\sum_{m=1}^{n-1} m\Delta^2 v(n-m+1).
    %
\end{align}


\section{Proofs for Weibull sources}
\label{secWeilbullA}

Proposition \ref{theoWeibull} seem to be new.

\subsection{Proof of Proposition \ref{theoWeibull}}

Let us introduce notations
\begin{align}
  \rho_w(t)&:=P\okra{N^W_w((0,t])>0},
  \\
  \rho_w(t|k)&:=P\okra{N^W_w((0,t])=k}.
\end{align}
We have $\mean V_W(t)=\sum_{w\in\mathbb{W}} \rho_w(t)$ and
$\mean V_W(t|k)=\sum_{w\in\mathbb{W}} \rho_w(t|k)$.

We choose a particular type $w$ and we suppress subscript $w$ in the
following.  We observe that
\begin{align}
  \rho(t|0)&=1-\rho(t)=e^{-(t/\lambda)^{\alpha}},
  \\
  \rho(t|k+1)&=\int_0^t \rho(t-s|k)\frac{d \rho(s)}{d s}d s.
\end{align}
For $\alpha\le 1$, we have $x^\alpha+y^\alpha\ge (x+y)^\alpha$. Hence
\begin{align}
  \rho(t|1)
  &=
    \int_0^t 
    \frac{\alpha}{\lambda}\okra{\frac{s}{\lambda}}^{\alpha-1}
    e^{-(s/\lambda)^{\alpha}-((t-s)/\lambda)^{\alpha}} d s
    \nonumber\\
  &\le
    e^{-(t/\lambda)^{\alpha}}
    \int_0^t 
    \frac{\alpha}{\lambda}\okra{\frac{s}{\lambda}}^{\alpha-1} d s
    \nonumber\\
  &=
    \okra{\frac{t}{\lambda}}^{\alpha}e^{-(t/\lambda)^{\alpha}}
    =
    \frac{t}{\alpha}\frac{d \rho(t)}{dt}.
\end{align}

In general, we have
\begin{align}
  \rho(t|k)\le \frac{\Gamma(\alpha+1)^k}{\Gamma(k\alpha+1)}
  \okra{\frac{t}{\lambda}}^{k\alpha}e^{-(t/\lambda)^{\alpha}}
\end{align}
since this pattern propagates by induction
\begin{align}
  \rho(t|k+1)
  &\le
    \frac{\Gamma(\alpha+1)^k}{\Gamma(k\alpha+1)}
    \int_0^t 
    \frac{\alpha}{\lambda}\okra{\frac{s}{\lambda}}^{\alpha-1}
    \okra{\frac{t-s}{\lambda}}^{k\alpha}
    e^{-(s/\lambda)^{\alpha}-((t-s)/\lambda)^{\alpha}} d s
    \nonumber\\
  &\le
    \frac{\Gamma(\alpha+1)^k}{\Gamma(k\alpha+1)}
    e^{-(t/\lambda)^{\alpha}}
    \int_0^t 
    \frac{\alpha}{\lambda}\okra{\frac{s}{\lambda}}^{\alpha-1}
    \okra{\frac{t-s}{\lambda}}^{k\alpha} d s
    \nonumber\\
  &=
    \frac{\Gamma(\alpha+1)^k}{\Gamma(k\alpha+1)}
    \okra{\frac{t}{\lambda}}^{(k+1)\alpha}
    e^{-(t/\lambda)^{\alpha}}
    \int_0^1 \alpha u^{\alpha-1}(1-u)^{k\alpha} d u
    \nonumber\\
  &=
    \frac{\Gamma(\alpha+1)^k}{\Gamma(k\alpha+1)}
    \frac{\Gamma(\alpha+1)\Gamma(k\alpha+1)}{\Gamma((k+1)\alpha+1)}
    \okra{\frac{t}{\lambda}}^{(k+1)\alpha}e^{-(t/\lambda)^{\alpha}}
    \nonumber\\
  &=
    \frac{\Gamma(\alpha+1)^{k+1}}{\Gamma((k+1)\alpha+1)}
    \okra{\frac{t}{\lambda}}^{(k+1)\alpha}e^{-(t/\lambda)^{\alpha}}.
\end{align}

Similarly, we can prove that for $\alpha\le 1$, we have
\begin{align}
  (-t)^k\frac{d^k \rho(t|0)}{dt^k}
  &=
    \sum_{j=0}^k w_{kj} \okra{\frac{t}{\lambda}}^{j\alpha} e^{-(t/\lambda)^{\alpha}},
\end{align}
where $w_{kj}\ge 0$ and $w_{kk}=\alpha^k$,
since this pattern propagates by induction
\begin{align}
  (-t)^{k+1}\frac{d^k \rho(t|0)}{dt^k}
  &=
  -t^{k+1}\frac{d}{dt}
    \frac{\sum_{j=0}^k w_{kj} \okra{\frac{t}{\lambda}}^{j\alpha} e^{-(t/\lambda)^{\alpha}}}{t^k}
    \nonumber\\
  &=
  \okra{k
  -
  t\frac{d}{dt}}
  \sum_{j=0}^k w_{kj} \okra{\frac{t}{\lambda}}^{j\alpha} e^{-(t/\lambda)^{\alpha}}
    \nonumber\\
  &=
    \sum_{j=0}^k \okra{k-j\alpha+\alpha\okra{\frac{t}{\lambda}}^{\alpha}}
    w_{kj} \okra{\frac{t}{\lambda}}^{j\alpha} e^{-(t/\lambda)^{\alpha}}.
\end{align}

Consequently, the claims follow by summing contributions $\rho(t)$ and
$\rho(t|k)$ over the types.

\renewcommand{\bibsep}{1pt}
\small
\bibliography{0-publishers-full,0-journals-full,books,ql,nlp,ai,math,mine}

\end{document}